\documentclass[11pt]{article}
\usepackage{amsmath,amssymb,latexsym}

\newcommand{\N}{\mathbb{N}}
\newcommand{\R}{\mathbb{R}}        

\newcommand{\limind}{\text{lim ind}}
\newcommand{\eps}{\varepsilon}

\newtheorem{theorem}{Theorem}
\newtheorem{lemma}[theorem]{Lemma}

\newcommand{\qed}{\hspace*{\fill} $\Box $}

\begin{document}
\title{A Fundamental System of Seminorms for $A(K)$}
\author{Dietmar Vogt}
\date{}
\maketitle

Let $K\subset\R^d$ be compact and $A(K)$ the space of germs of real analytic functions on $K$ with its natural (LF)-topology (see e.g. \cite{MV}, 24.38, (2)). This topology can also be given by $A(K)=\limind_{k\to+\infty} A_k$ where
$$A_k=\{(f_\alpha)_{\alpha\in\N_0^d}\in C(K)^{\N_0^d}\,:\, \|f\|_k:=\sup_{x\in K} \frac{|f^{(\alpha)}(x)|}{\alpha!} k^{-|\alpha|}< +\infty\}.$$
Based on this description we give in the present note an explicit fundamental system of seminorms for $A(K)$.

We start with a modified problem.

Let $X$ be a Banach space. We put
$$F_k=\{(x_\alpha)_{\alpha\in\N_0^d}\in X^{\N_0^d}\,:\, \|x\|_k:=\sup_\alpha \|x_\alpha\| k^{-|\alpha|}< +\infty\}$$
and
$$F=\limind_{k\to+\infty} F_k.$$
On $F$ we consider for any positive null-sequence $\delta=(\delta_n)_{n\in\N}$ the continuous norm
$$|x|_\delta = \sup_\alpha \|x_\alpha\| \delta_{|\alpha|}^{|\alpha|}.$$
\begin{lemma} \label{lem}The norms $|\cdot|_\delta$ are a fundamental system of seminorms on $F$.
\end{lemma}

\bf Proof: \rm It is sufficient to show  that for every positive sequence $\eps_k$, $k\in\N$, there exists $\delta$ such that
$$U_\delta:=\{x\,:\,|x|_\delta\le 1\}\subset\sum_k \eps_k B_k$$
where $B_k$ denotes the unit ball of $F_k$. Without restriction of generality we may assume, that $\eps_k\le1$ for all $k$.

Fot every $k$ we choose $n_k>n_{k-1}$, such that $k-1<\eps_k^{1/n_k}k$. We put $\delta_n^{-1} = \eps_k^{1/n_k}k$ for $n_k\le n< n_{k+1}$. We obtain for these $n$
$$\delta_n^{-n}=\eps_k^{n/n_k}k^n\le \eps_k k^n.$$
Due to the construction $\delta=(\delta_n)_n$ is a null-sequence. For $x\in U_\delta$ and  $n_k\le |\alpha|<n_{k+1}$ we have $x_\alpha\in \eps_k B_k$  and therefore
$$\xi_k=\sum_{n_k\le|\alpha|<n_{k+1}} x_\alpha\in \eps_k B_k.$$
Since $x=\sum_k \xi_k$ the proof is complete. \qed

\bigskip

\begin{theorem} If $K\subset \R^d$ is compact, then the norms
 $$|f|_\delta = \sup_\alpha \sup_{x\in K} \frac{|f^{(\alpha)}(x)|}{\alpha!} \delta_{|\alpha|}^{|\alpha|},$$
where $\delta$ runs through all positive null-sequences, are a fundamental system of seminorms in $A(K)$.
\end{theorem}

\bf Proof: \rm Let $F$ be as above with $X=C(K)$. We define a map $A:A(K)\to F$ by $A(f)=\Big(\frac{f^{(\alpha)}(x)}{\alpha!}\Big)_{\alpha\in\N_0^d}$. The map $A$ is obviously continuous and $A^{-1}(B)$ is bounded in $A(K)$ for every bounded set $B$ in $F$. From Baernstein's Lemma (see \cite{MV}, 26.26) it follows, that $A$ is an injective topological homomorphism. Hence Lemma \ref{lem} proves the result. \qed

\vspace{.5cm}

\noindent Bergische Universit\"{a}t Wuppertal,
\newline FB Math.-Nat., Gau\ss -Str. 20,
\newline D-42119 Wuppertal, Germany
\newline e-mail: dvogt@math.uni-wuppertal.de

\end{document}